%
%
%
%
%
%
%
\input amstex
%
%

\def\A{\Cal{A}}
\def\B{\Cal{B}}
\def\C{\Cal{C}}

\def\H{\Cal{H}}

\def\T{\Cal{T}}

\def\N{\Cal{N}}
\def\to{\rightarrow}
\def\Prim{\hbox{Prim}}

\def\hull{\hbox{hull}}
\def\ker{\hbox{ker}}
\def\qedd{
\hfill
\vrule height4pt width3pt depth2pt
\vskip .5cm
}
\NoBlackBoxes
%
\documentstyle{amsppt}
%
\loadbold
\magnification=\magstep 1
%
%
\topmatter
\title The topology of ideals in some triangular AF algebras
\endtitle
\author Michael P.\ Lamoureux
\endauthor
\affil University of Calgary
\endaffil
\address
Department of Mathematics and Statistics,
University of Calgary,
2500~University Drive N.W.,
Calgary~AB, Canada T2N 1N4.
\endaddress
\email mikel\@ math.ucalgary.ca
\endemail
\date December 1995
\enddate
\keywords Almost finite algebra, lattice of ideals, meet-irreducible, nests of
subspaces
\endkeywords
\subjclass
 46K50, 47D25, 46M40
\endsubjclass

\abstract
%
%
A  set of meet-irreducible ideals is described for a class of maximal
triangular almost
finite algebras. This set forms a topological space under the hull-kernel
closure, and there is a
one-to-one correspondence between closed sets in this space and ideals in the
AF algebra. Some
connections with nest representations and nest-primitive ideals are also
described.
\endabstract

\endtopmatter
\document


\head Introduction
\endhead

There are many instances in the study of operator algebras where one can
describe the ideal
structure of an algebra by identifying the set of ideals with a collection of
closed sets in some
topological space. The canonical example is in C*-theory, where one may
construct the primitive
spectrum $\Prim(\A)$ of a C*-algebra $\A$ as the set of primitive ideals in
$\A$, and
introduce the the Jacobson (or hull-kernel) topology to create a topological
space whose closed
sets are in one-to-one correspondence with the closed ideals in $\A$. More
precisely, in this
case every ideal is described as the intersection of all those primitive ideals
which contain
it. When the algebra $\A$ arises as the C*-crossed product of a group acting on
a space
$X$, there is a close connection between the topology of $\Prim(\A)$ and the
structure of the
orbits on $X$, and so the ideal structure can be understood in terms of these
orbits. (See
for instance [4] or [17].) In [7] and [8], this notion is extended to
nonselfadjoint crossed
products, where it is shown that  so-called nest-primitive ideals corresponding
to arcs
(subintervals of orbits) in
$X$ could be used to build a topological space which again describes the total
ideal structure of
the nonselfadjoint algebra.  Indeed, in this case it was the
structures on $X$ and their topology which provided the motivation for studying
this class of
ideals.

However, it seems worthwhile to ask directly the question: given an algebra of
operators, what
set of ideals can be used to build a topological space, and when is this space
rich enough to
describe all the ideals in the algebra? Having asked the question, we will show
via
Propositions~1.1 and 1.3 that the ideals of interest are precisely the
meet-irreducible ideals;
that is, those ideals which cannot be written as the intersection of two
strictly larger ideals.
For C*-algebras, these are nothing more than the primitive ideals (see
Proposition~2.1), so the
more interesting examples occur for nonselfadjoint operator algebras.

One area where much work has been done to understand the ideal structure for
nonselfadjoint
algebras is in the study of
direct limits of upper triangular matrix algebras, and more generally
triangular almost finite
algebras, where questions about ideal
structure have been an important consideration in understanding the
classification of the
algebras. For instance, in [14] Power  has identified both the lattice of all
ideals as well as the
meet-irreducible ideals for an infinite tensor product of upper triangular
matrix algebras. As
noted in [14], this lattice of ideals is distributive, thus by
Propositions~1.1 and 1.3 of the present work,  the set of meet-irreducible
ideals is in fact a
topological space. Moreover, the closed sets in this space are in one-to-one
correspondence with
closed ideals in the tensor product, for every ideal is the intersection of the
meet-irreducible
ideals which contain it, as can be readily verify from the constructions in
[14]. Other authors
(as in [3],[6] and [12]) have looked at the lattice structure of ideals in
other AF algebras,
including the dual notion of join-irreducibility (for instance, in [6]); one
might expect from the
C*-case that such ideals are rather rare. Meet-irreducible ideals also appeared
in [7] and [8], in
the context of nest-primitive ideals in nonselfadjoint crossed products, but in
a rather curious
way: the ideals constructed were all nest-primitive, and happened to form a
topological space; by
Propositions~1.1 and 1.3, these are necessarily meet-irreducible.

The purpose of this paper is to describe explicitly how a topology arises on a
set of ideals, and
then apply this to a class of maximally triangular AF algebras. The paper is
structured as
follows. Section One identifies the conditions necessary on a set of ideals to
become a
topological space. Section Two considers a number of motivating examples;  for
C*-algebras, we see the meet-irreducible ideals are just the primitive ideals,
while for two
types of nonselfadjoint algebras (the compact nest algebra, and the disk
algebra), the
meet-irreducibles are a generalization of the primitive ideals, namely the
nest-primitive ideals.
Section Three contains the core of the paper, where a large set of
meet-irreducible ideals is
identified in a strongly maximal triangular AF algebra, which produces a
topological space whose
closed sets are in one-to-one correspondence with ideals in the AF algebra.

Throughout the paper, we will discuss the connection of these meet-irreducible
ideals with the
nest-primitive ideals, which were introduced in [7] to describe the topology of
the
nonselfadjoint crossed products. However, except for some motivating examples,
it is not clear what
general relationship exists between these two types of ideals.

\vskip .3cm

Our notation is as follows. An almost finite (AF) C*-algebra $\B$ is the
norm-closed union of an
increasing sequence of finite dimensional C*-algebras $\B_0\subseteq
\B_1\subseteq\ldots$, and a
triangular AF algebra is a norm-closed subalgebra $\A \subseteq \B$ with
$\A\cap\A^*$ a canonical
masa in $\B$. (That is, with the diagonal $\C = \A\cap\A^*$, we have each $\C_k
= \C\cap\B_k$  a
masa in $\B_k$, and $\C$ is the norm closure of the union of the $\C_k$.) We
say $\A$ is strongly
maximal triangular if each finite algebra
$\A_k =
\A\cap\B_k$ is maximal triangular in
$\B_k$; equivalently, we may assume the algebras $\B_k$ are chosen so that
$\A_k \subseteq \B_k$
is a direct sum of upper triangular matrix algebras. (See [14] for a
well-developed account of
these algebras.)

We will make use of the fact that the set of closed, two-sided ideals in a
Banach algebra form a
multiplicative lattice under the inclusion ordering; the greatest lower bound
(meet) of two
ideals is their intersection, the least upper bound (join) is the closed linear
span of the union,
and the product is the closed linear span of pairwise element products. A
maximal ideal is thus a
proper ideal which is maximal in the inclusion order, while a primary ideal is
a proper ideal
which is contained in a unique maximal ideal. Throughout this paper, an ideal
means a norm-closed,
two-sided ideal.

A bounded representation $\pi$ of a Banach algebra $\A$ on Hilbert space $\H$
is called a {\it nest
representation} if the closed, $\pi(A)$-invariant subspaces in $\H$ are
linearly ordered by
inclusion; that is, they form a nest (see [7]). This is a generalization of the
notion of an
irreducible representation, where the invariant subspaces form the trivial nest
$\{ 0, \H \}$. The
kernel of an irreducible representation is a primitive ideal; the kernel of a
nest representation
will be called a {\it nest-primitive} (or n-primitive) ideal.

\head I. Topological spaces of ideals
\endhead

Given a set $\Omega$ of closed, two-sided ideals in a Banach algebras $A$, we
can attempt to
introduce a topology on $\Omega$ by defining formally the closure of a subset
$F$ in
$\Omega$ as
$$ \overline{F} = \hull(\ker (F))$$
where $I = \ker (F)$ is the ideal obtained as the intersection of all the
ideals in set $F$ and
$\hull (I)$ is the subset of ideals in $\Omega$ which contain $I$. To ensure
that this operation
gives a topology, it suffices to verify that the four Kuratowski axioms for a
closure operation
are satisfied. Namely, the closure of the empty set $\phi$ is empty,
\roster
\item"{(K1)}" $\overline{\phi} = \phi$;
\endroster
and for any subsets $F,G$ in $\Omega$, we have
\roster
\item"{(K2)}" $\overline{F} \supseteq F$;
\item"{(K3)}" $\overline{\overline{F}} = \overline{F}$;
\item"{(K4)}" $\overline{F\cup G} = \overline{F}\cup\overline{G}$.
\endroster
It is not hard to see that $\overline{\phi} = \phi$ if and only if $\Omega$
consists of proper
ideals. The equality $\overline{\overline{F}} = \overline{F}$ and the
containments
$\overline{F} \supseteq F$ and $\overline{F\cup G} \supseteq
\overline{F}\cup\overline{G}$ follow
almost immediately from definition of hull and ker, so to obtain a topology,
one only requires
the reverse containment for (K4). We state a convenient equivalency in the
following proposition.

\proclaim{Proposition 1.1}
Let $\Omega$ be a set of proper ideals in algebra $\A$. Then the hull-kernel
operation produces a
topological closure operation if and only if for any ideal $I$ in $\Omega$ and
subsets $F,G$ we
have
$$I\supseteq \ker (F)\cap\ker (G) \quad\hbox{ implies }\quad
I\supseteq \ker (F) \hbox{ or } I\supseteq \ker (G).$$
\endproclaim
\demo{Proof}
By the comments above, we obtain a topological closure if and only if we obtain
the
containment in (K4) of $\overline{F\cup G} \subseteq
\overline{F}\cup\overline{G}$; that is
$$ I \in \overline{F\cup G} \Rightarrow I\in \overline{F}\cup\overline{G}.$$
But by definition, $I \in \overline{F\cup G}$ is equivalent to $I\supseteq \ker
(F) \cap \ker (G)$
while
$I\in \overline{F}\cup\overline{G}$ is equivalent to $I\supseteq \ker (F)$ or
$I\supseteq \ker
(G)$. Thus (K4) holds if and only if
$$I \supseteq \ker (F) \cap \ker (G) \Rightarrow I\supseteq \ker (F) \hbox{ or
}I\supseteq \ker
(G).$$
\qedd
\enddemo

As our desire is to find a set $\Omega$ of ideals which is rich enough that
every ideal in the
algebra can be written as $\ker(F)$ for some subset $F$ of $\Omega$, it is not
too restrictive to
generalize the condition in the above proposition as follows.

\proclaim{Definition 1.2}
An ideal $I$ in algebra $A$ is called (K4) if, for any ideals $J,K$ in the
algebra, one has that
$$I\supseteq J\cap K \Rightarrow I\supseteq J \hbox{ or } I\supseteq K.$$
\endproclaim

Note that the condition (K4) is similar to the condition that an ideal be prime
(i.e. $I\supseteq J\cdot K \Rightarrow I\supseteq J \hbox{ or } I\supseteq K$)
as well as the
condition that it be meet-irreducible
(i.e. $I=J\cap K \Rightarrow I = J \hbox{ or } I = K$). We use the terminology
(K4) because of
the connection with Kurotowski's fourth axiom; however, in light of the
following proposition,
the term ``meet-prime'' might be appropriate as well.

\proclaim{Proposition 1.3} For any closed, two-sided ideal $I$  in
Banach algebra $A$,
\roster
\item if $I$ is prime, then  $I$ is (K4);
\item if $I$ is (K4), then $I$ is meet-irreducible.
\endroster
\endproclaim
\demo{Proof}
If $I$ is prime and $I\supseteq J\cap K$, then $I\supseteq J\cap K \supseteq
J\cdot K$, so by
primeness, either
$I\supseteq J$ or $I\supseteq K$, hence $I$ is (K4).

If $I$ is (K4) and $I = J\cap K$, then $J,K \supseteq I\supseteq J\cap K$, and
so
either
$J \supseteq I\supseteq J$ or $K \supseteq I\supseteq K$. That is, $I=J$ or
$I=K$, so $I$ is
meet-irreducible.
\qedd
\enddemo

Note that the conditions prime and (K4) are quite distinct. For instance, in
the algebra of upper
triangular four by four matrices, the ideal
$$ I = \left( \matrix
 * & * & * & *\\
   & 0 & 0 & *\\
   &   & 0 & *\\
   &   &   & *
\endmatrix \right)$$
is (K4) but not prime.
On the other hand, meet-irreducible and (K4) are equivalent notions for
distributive lattices of
ideals, for if $I$ is meet-irreducible, and $I\supseteq J\cap K$, then using
$\vee$ for the
joint operation, $I = (J\cap K)\vee I = (J\vee I)\cap(K\vee I)$, so by
irreducibility, $I = J\vee
I$ or
$I = K\vee I$, which means $I$ contains $J$ or $I$ contains $K$. That is, $I$
is (K4).

However, if the lattice of ideals can be described by closed sets in a
topological space, then
the lattice is distributive (since the closed sets in a topological space form
a distributive
lattice, under the operations of union and intersection), so if we have any
hope of finding a
``good'' space of ideals, we need to look at the meet-irreducibles.

We will make use of the following connections between an ideal and its quotient
algebra. Note this
is a standard result when applied to prime ideals in C*-algebras.

\proclaim{Proposition 1.4} Let $I$ be a closed, two-sided ideal in Banach
algebra $A$.
\roster
\item If $I$ is prime in $A$, then $0$ is prime in $A/I$;
\item if $I$ is (K4) in $A$, then $0$ is (K4) in $A/I$;
\item if $I$ is meet-irreducible in $A$, then $0$ is meet-irreducible in $A/I$.
\endroster
\endproclaim
\demo{Proof}
The three cases are demonstrated via the same technique. We will show (3).
Assume $I$ is meet-irreducible in $A$. If $0$ is not meet-irreducible in $A/I$,
then there exist
two non-zero ideals $J_0$ and $K_0$ in $A/I$ with zero intersection. Thus,
there are two ideals
$J$ and $K$ in $A$ properly containing $I$ with $J/I = J_0$ and $K/I = K_0$.
Notice
$(J\cap K)/I = J_0\cap K_0 = 0$ and so $J\cap K = I$, contradicting that $I$
was
meet-irreducible.
\qedd
\enddemo

The following section describes several examples where the connections between
ideals types can
be made quite clear, and are the motivation to consider more general triangular
algebras.

\head II. C*-algebras and compact nest algebras
\endhead

For C*-algebras, there is an equivalence of five conditions.

\proclaim{Theorem 2.1} Let $I$ be a closed, two-sided ideal in a separable
C*-algebra. Then the
following are equivalent:
\roster
\item $I$ is nest-primitive;
\item $I$ is primitive;
\item $I$ is prime;
\item $I$ is (K4);
\item $I$ is meet-irreducible.
\endroster
Moreover, every ideal in a separable C*-algebra is equal to the intersection of
meet-irreducible
ideals.
\endproclaim
\demo{Proof}
The equivalence of (1) and (2) is shown in the comments before Theorem~I.4 of
[7];
to summarize, since a nest representation is cyclic, one can use a result of
U.\ Haagerup to
conclude it is similar to a *-representation, and thus the nest of invariant
subspaces is
trivial. In particular, a nest representation is similar to an irreducible
*-representation, and
thus its kernel is a primitive ideal.

The implication  (2) $\Rightarrow$ (3) is well-known for C*-algebras; see for
instance
Prop.~3.13.10 of [9]. The converse is true for separable C*-algebras, as shown
in
Prop.~4.3.6 of the same text.

The implications (3) $\Rightarrow$ (4) $\Rightarrow$ (5) follow from
Proposition~1.3. To complete
the proof, we will show that (5) implies (2).

To this end, suppose the ideal $I$ is meet-irreducible. The set $\Prim(A/I)$ of
primitive ideals
in the quotient C*-algebra $A/I$ forms a topological space in the hull-kernel
topology, which is
second countable since
$A/I$ is separable. Let $\{ G_n \}$ be a countable basis for the topology,
consisting of open,
non-empty sets in $\Prim(A/I)$.  As
$G_n$ is non-empty, the ideal $\ker(\Prim((A/I)\backslash G_n))$ is non-zero,
while the
intersection
$\ker(G_n)\cap\ker(\Prim((A/I)\backslash G_n))$ is the zero ideal,  which is
meet-irrreducible by
Proposition~1.4. Thus $\ker(G_n)$ is the zero ideal, so each set $G_n$ is dense
in the space
$\Prim(A/I)$, which is a Baire space. By the Baire category theorem, $\cap G_n$
is non-empty and
contains a point $J$ which is dense in $\Prim(A/I)$. Hence $J$ is an ideal
contained in every
primitive ideal, so in fact $J = 0$, consequently $0$ is a primitive ideal.
Letting $\pi$ be an
irreducible representation of $A/I$ with 0 kernel, its pullback to $A$ is an
irreducible
representation with kernel $I$. Hence $I$ is primitive.
\qedd
\enddemo

\demo{Remark}
It is interesting to note that the proof of (5) implies (2) above uses the same
techniques as in
the demonstration of the well-known result in C*-theory that  prime and
primitive are
equivalent for separable C*-algebras. Moreover, by the correspondence between
closed ideals in a C*-algebra $A$ and closed sets in the primitive ideal space,
it is clear that
the lattice structure is isomorphic to the lattice of closed sets in
$\Prim(A)$. Meet-irreducible ideals correspond to minimal closed sets, namely
single points.
Join-irreducible ideals correspond to maximal closed sets, namely those that
cannot be expressed
as the intersection of two larger closed sets. In the Hausdorff case, this
corresponds to the set
compliment of an isolated point.
In particular, while there are lots of meet-irreducible ideals around,
join-irreducible ideals
are rather rare.
\enddemo

We now consider two nonselfadjoint examples, where the description of the
meet-irreducible ideals
can also be made precise.
The interest in these compact nest algebras and the disk algebra rises from
their connection with
nonselfadjoint crossed products of dynamical systems. Their are, in some sense,
the elementary
building blocks for the crossed products, as described in [7] and [8].

Recall a nest in a Hilbert space $\H$ is a family of closed subspaces
which is linearly ordered by inclusion, contains $0$ and $\H$, and is closed
under arbitrary
unions and intersections. An interval is the orthogonal difference $N_1 \ominus
N_0$ of two
elements of the nest. A nest algebra is the set of all bounded linear operators
on
$\H$ leaving invariant the elements of a fixed nest. We say a compact nest
algebra is the set of
all compact linear operator on $\H$ leaving invariant a fixed nest.
Notice the identity representation of a compact nest algebra $\A$ on a Hilbert
space is a nest
representation, as is it compression to any interval $N_1\ominus N_0$, for the
invariant subspaces
are just the nest itself, or its restriction to the interval. (See Theorems~I.4
and II.1 of [7]).
The nest representations and
nest-primitive ideals of an arbitrary compact nest algebra can thus be
described completed, and we
observe the close relations between nest-primitive ideals and
meet-irreducibles.

\proclaim{Proposition 2.2} Let $\A$ be a compact nest algebra on Hilbert space
$\H$, and $I$ a
closed, two-sided ideal in $\A$. Then the following are equivalent:
\roster
\item $I$ is nest-primitive;
\item $I$ is (K4);
\item $I$ is meet-irreducible;
\item $I$ is kernel of the compression map of $\A$ to some interval of the
nest.
\endroster
Moreover, every ideal in $\A$ is an intersection of the meet-irreducible ideals
which contain it.
\endproclaim
\demo{Proof}
These results are  a summary of work in [7]. The equivalence
$(4)\Leftrightarrow
(1)$ and the implication $(1)\Rightarrow (2)$ follow from  Theorems~II.1 and
II.2 in [7],
while $(2)\Rightarrow(3)$ is Corollary~II.3, although this is also a general
consequence of Proposition~1.3 of the present paper. The key construction in
these results is the
observation that every ideal in $\A$ is of the form
$$I_\phi = \{ a\in \A : aN \leq \phi(N), \hbox{ for all } N \in \N \},
$$
where $\phi:\N \rightarrow\N$ is a left continuous order homomorphism of the
nest $\N$ into
itself, with $\phi(N)\leq N$ for every subspace $N$ in the nest. The
nest-primitive ideals are
determined by a pair of subspaces $N_0\leq N_1$, with corresponding order
homomorphism
$$\phi(N) = \cases N,&N<N_0\\ N_0, &N_0\leq N\leq N_1\\ N, &N>N_1 .\endcases
$$
That is, $\phi$ determines a corner, and $I_\phi$ is the kernel of the
compression of the
identity representation onto the interval $N_1\ominus N_0$ in the nest.

To complete the proof of the proposition,, we need to show that (3) implies
(4). Assume that the
ideal
$I$ is meet-irreducible. If $I = I_\phi$ is not the kernel of a compression
map, then $\phi$ is
not of the form of a corner, as above, so we may construct left continuous
order homomorphisms
$\phi_1,\phi_2$ not equal to $\phi$ with
$$\phi = \max\{\phi_1,\phi_2\}.
$$
Thus $I_\phi = I_{\phi_1}\cap I_{\phi_2}$, contradicting that $I$ is not
meet-irreducible. Hence
$I$ is indeed the kernel of a compression.

The fact that every ideal in $\A$ is the intersection of meet-irreducibles
(equivalently,
nest-primitives) is Corollary~II.4 of [7].
\qedd
\enddemo

We have immediately the following.

\proclaim{Corollary 2.3} Let $\A$ be a compact nest algebra and $\Omega$ the
set of
meet-irreducible ideals in $\A$. Then the hull-kernel operation gives a
topology on $\Omega$, and
there is a one-to-one correspondence between ideals in $\A$ and closed sets in
$\Omega$, given by
the hull and kernel maps.
\endproclaim

The paper [7] also considers the disk algebra $A(\Delta)$, the Banach algebra
of
continuous function on the unit disk in the complex plane, which are analytic
in the interior.
This is an interesting example where the meet-irreducible ideals are not a rich
enough set to
describe all the ideals in the algebra. Indeed, it is well-known that  every
ideal in
$A(\Delta)$ is principal and is generated by some inner function in
$A(\Delta)$, with a singular
part determined by a measure on the unit circle. Maximal ideals correspond to
ideals of functions
vanishing at a single point in the disk; thus primary ideals have a singleton
for their zero set.
However,  not every ideal is the intersection of primary ideals; we miss the
ones with
non-discrete measures in the singular part of the generating function (see
[5]). We may
extend the relevant results in [7] as follows.

\proclaim{Proposition 2.4} Let $A(\Delta)$ be the disk algebra, and $I$ a
closed, two-sided ideal
in $A(\Delta)$. Then the following are equivalent:
\roster
\item $I$ is nest-primitive;
\item $I$ is (K4);
\item $I$ is meet-irreducible;
\item $I$ is either primary, or zero.
\endroster
However, not every ideal is the intersection of meet-irreducible ideals.
\endproclaim
\demo{Remark} This proposition is essentially Theorems~IV.1 and IV.2 of [7],
although
these only show
$(1)\Leftrightarrow(4)\Rightarrow(2)$, while the implication
$(2)\Rightarrow(3)$ follows from
Proposition~1.3. To see that $(3)\Rightarrow(4)$, the same argument as used in
Theorem~IV.1 of
[7]. That is, if $I$ is a non-zero meet-irreducible ideal, then $I$ is
generated by some
function $F$ in $A(\Delta)$. The zero set of $F$ must be a singleton (and hence
$I$ is primary),
for if not, we factor $F = fg$ as a product of two functions with strictly
smaller zero sets.
Then if $J$ denotes the ideal generated by $f$ and $K$ the ideal generated by
$g$, we see that
$fg \in J\cap K$, so $I = J\cap K$, contradicting that $I$ is meet-irreducible.
\enddemo

\head III. Finite and almost finite triangular algebras
\endhead

The algebra $\T_n$ of $n$ by $n$ upper triangular matrices is a special case of
a compact nest
algebra and thus is covered by Proposition~2.2 above. Note that the typical
nest-primitive ideal is
formed by a corner of zeroes, such as
$$ I = \left( \matrix
 * & * & * & *\\
   & 0 & 0 & *\\
   &   & 0 & *\\
   &   &   & *
\endmatrix \right),$$
and thus the nest-primitive ideals are in one-to-one correspondence with matrix
units in $\T_n$.
To make this correspondence explicit, fix a system of matrix units $\{ e_{ij}
\}_{i,j=1}^n$ for
$\T_n$, and for any unit
$e_{ij}$ from the system, let $I(e_{ij})$ denote the largest ideal in $\T_n$
not containing
$e_{ij}$. Conversely, given a meet-primitive ideal $I$ in $\T_n$, let $e(I)$
denote the
matrix unit in the system of units which is not contained in $I$, and maximal
of all such
units with respect to the partial order
$$ e_{ij} \leq_p e_{i'j'}  \quad\Leftrightarrow\quad  j\leq j'  \hbox{ and }
i\geq i'.
$$
This order $\leq_p$ is closely related to the Peters-Poon-Wagner order  on the
diagonal (see [12]),
where two projections $p,q$ are ordered $p \preceq q$ if there is some element
$w$ in the
normalizer of the diagonal with range projection $r(w) = p$ and domain
projection $d(w) = q$. Thus
we may write for two units $e,f$ in $\T_n$ that
$$ e \leq_p f  \quad\Leftrightarrow\quad d(e)\preceq d(f) \hbox{ and }
r(e)\succeq r(f).
$$
In terms of the nest-primitive ideals $I(e), I(f)$ defined above, we may also
write
$$ e \leq_p f  \quad\Leftrightarrow\quad  I(e) \supseteq I(f).
$$

This description extends in a straightforward manner to a finite direct sum of
upper triangulars.
Namely, if $\A = \T_{n_1} \oplus \cdots \oplus \T_{n_k}$ is a direct sum of
upper triangular
matrices of sizes $n_1, \ldots n_k$, we fix a system of matrix units $\Cal E$
and introduce the
partial order
$$ e \leq_p f  \quad\Leftrightarrow\quad  d(e)\preceq d(f) \hbox{ and }
r(e)\succeq r(f)
$$
for any two units $e,f$ in $\Cal E$. For such a unit $e$, let $I(e)$ denote the
largest ideal in
$\A$ not containing $e$. That is, $I(e)$ is the closed linear span, or join, of
all ideals in
$\A$ not containing $e$. Since every matrix unit $e$ can be taken in the form
$0\oplus\cdots\oplus e^j\oplus\cdots\oplus 0$ in
$\T_{n_1}\oplus\cdots\oplus\T_{n_k}$, where $e^j$
is a matrix unit in
$\T_{n_j}$, it is easy to see that
$$I(e) =  \T_{n_1}\oplus\cdots\oplus I(e^j) \oplus\cdots\oplus\T_{n_k},$$
where $I(e^j)$ is the meet-irreducible ideal in $\T_{n_j}$ corresponding to the
unit $e^j$.

\proclaim{Proposition 3.1}
The ideal $I(e) = \T_{n_1}\oplus\cdots\oplus I(e^j) \oplus\cdots\oplus\T_{n_k}$
is
meet-irreducible in $\A = \T_{n_1} \oplus \cdots \oplus \T_{n_k}$, and every
meet-irreducible ideal
is of this form.
\endproclaim
\demo{Proof}
We will show $I(e)$ is (K4) and hence meet-irreducible by Proposition 1.3.
Suppose $I(e)\supseteq
J\cap K$, where $J,K$ are ideals in the algebra $\A$. Since $e\not\in I(e)$, we
must have either
$e\not\in J$ or $e\not\in K$. Since $I(e)$ is the join of all ideals not
containing $e$, we have
either $I(e) \supseteq J$ or $I(e) \supseteq K$. Thus $I(e)$ is (K4) and hence
meet-irreducible.

Conversely, suppose $I$ is a meet-irreducible ideal in $\A$. Since $\A$ is a
direct sum, we may
write
$I = I_1\oplus\cdots\oplus I_k$, where each $I_j$ is an ideal in the
corresponding $\T_{n_j}$. If
two or more of the $I_j$ are proper ideals, say $I_{j_1} \neq \T_{n_{j_1}}$ and
$I_{j_2} \neq
\T_{n_{j_2}}$, then we can write
$$\align
 I = &\quad (I_1\oplus\cdots\oplus\T_{n_{j_1}}\oplus\cdots\oplus
I_{j_2}\oplus\cdots\oplus I_k) \cr
     &\cap (I_1\oplus\cdots\oplus\
I_{j_1}\oplus\cdots\oplus\T_{n_{j_2}}\oplus\cdots\oplus I_k) \cr
\endalign
$$
contradicting that $I$ is meet-irreducible. Thus all but at most one of the
ideals $I_j$ equal the
full algebra $T_{n_j}$, so we write
$$ I = \T_{n_1}\oplus\cdots\oplus I_j\oplus\cdots\oplus \T_{n_k},
$$
where $I_j$ is a (possibly proper) ideal in $\T_{n_j}$.

If $I_j$ is not meet-irreducible in $T_{n_j}$, we can write $I_j = J_j\cap K_j$
for two strictly
larger ideals $J_j$ and $K_j$ in $T_{n_j}$, so
$$ I = (\T_{n_1}\oplus\cdots\oplus J_j\oplus\cdots\oplus \T_{n_k})\cap
       (\T_{n_1}\oplus\cdots\oplus K_j\oplus\cdots\oplus \T_{n_k}),
$$
again contradicting that $I$ is meet-irreducible.

Thus $I_j$ is meet-irreducible in $\T_{n_j}$, so by Proposition~2.2, we have
$I_j = I(e^j)$ for
some matrix unit in
$\T_{n_j}$, and thus
$$\align I &= \T_{n_1}\oplus\cdots\oplus I(e^j) \oplus\cdots\oplus \T_{n_k} \cr
           &= I(0\oplus\cdots\oplus e^j \oplus\cdots\oplus 0) \cr
           &= I(e),
\endalign
$$
where $e = 0\oplus\cdots\oplus e^j \oplus\cdots\oplus 0$ is a matrix unit for
$\A$.
\qedd
\enddemo

Notice this proof also shows an ideal in $\T_{n_1}\oplus\cdots\oplus\T_{n_k}$
is (K4) if and only
if it is meet-irreducible, and is thus of the form $I(e)$, for some matrix unit
$e$. It is also
easy to see that $I(e)$ is the kernel of a nest representation, namely the one
obtained by
compressing the natural component representation
$$\T_{n_1}\oplus\cdots\oplus\T_{n_k} \rightarrow \T_{n_j} \rightarrow \B({\Bbb
C}^{n_j})$$
to an interval in ${\Bbb C}^{n_j}$, corresponding to the matrix unit $e^j$ in
$e =
0\oplus\cdots\oplus e^j\oplus\cdots\oplus 0$. Thus $I(e)$ is also a
nest-primitive ideal.

Conversely, if $I = I_1\oplus\cdots\oplus I_k$ is a nest-primitive ideal, then
it is the kernel
of some nest representation $\pi$. If two or more of the $I_j$ are proper,
say $I_{j_1} \neq \T_{n_{j_1}}$ and $I_{j_2} \neq \T_{n_{j_2}}$, then the
ideals
$J = 0\oplus\cdots\oplus\T_{n_{j_1}}\oplus\cdots 0$ and
$K  = 0\oplus\cdots\oplus\T_{n_{j_2}}\oplus\cdots 0$ are not in the kernel of
$\pi$. Since $\pi$ is
a nest representation, the non-zero invariant subspaces $\overline{\pi(J)\H}$
and
$\overline{\pi(K)\H}$ are ordered, so we may assume $0\neq \overline{\pi(J)\H}
\subseteq
\overline{\pi(K)\H}$.  Multiplying by $\pi(J)$, we see
$$ 0 \neq  \overline{\pi(J^2)\H} \subseteq \overline{\pi(JK)\H} = 0,
$$
a contradiction. Hence, at most only one of the $I_j$ can be proper. Thus
$I = T_{n_1}\oplus\cdots\oplus I_j \oplus\cdots\oplus T_{n_k}$, and applying
Proposition 2.2 to
the component $\T_{n_j}$, we conclude $I_j = I(e^j)$ for some matrix unit $e^j$
in $T_{n_j}$.
Thus,
$$\align I &= \T_{n_1}\oplus\cdots\oplus I(e^j) \oplus\cdots\oplus \T_{n_k} \cr
           &= I(0\oplus\cdots\oplus e^j \oplus\cdots\oplus 0) \cr
           &= I(e).
\endalign
$$
That is, the nest-primitive ideals are all of the special form constructed
above. We may
summarize as follows.

\proclaim{Proposition 3.2}
Let $\A = \T_{n_1}\oplus\cdots\oplus\T_{n_k}$ be a direct sum of upper
triangular $n_j$ by $n_j$
matrix algebras, and $I$ a two-sided ideal in $\A$. Then the following are
equivalent:
\roster
\item $I$ is nest-primitive;
\item $I$ is (K4);
\item $I$ is meet-irreducible;
\item There exists a matrix unit $e$ in $\A$ with
$I = I(e) = \overline{\hbox{span}}\{ J \hbox{ an ideal } : e\not\in J\}$.
\endroster
Moreover, for a given system of matrix units $\Cal E$ in $\A$, there is a
one-to-one
correspondence between nest-primitive ideals in $\A$ and matrix units in $\Cal
E$, given by (4).
\endproclaim

{}From this, we may now build a finite topological space which describes the
lattice of ideals for
$\A$.

\proclaim{Proposition 3.3}
Let $\A = \T_{n_1}\oplus\cdots\oplus\T_{n_k}$ be a direct sum of upper
triangular matrices, and
$\Omega$ the set of nest-primitive (equivalently, meet-irreducible) ideals in
$\A$. Then the
hull-kernel closure produces a topology on $\Omega$, and there is a one-to-one
correspondence
between closed sets
$F \subseteq \Omega$ and ideals $I \in \A$, given by
$$
\align                    I &= \ker(F)   \cr
       \hbox{ and }\quad  F &= \hull(I).
\endalign
$$
\endproclaim
\demo{Proof}
The fact that the meet-irreducible ideals are (K4) provides the topology, by
Proposition 1.1. We
only need to show that every ideal in $\A$ is the intersection of the
meet-irreducible ideals
which contain it. But, any ideal $J$ is determined by the matrix units it
contains, equivalently
by the units it does not contain, so
$$\align
J &= \cap\{ I(e) : e\not\in J \} \cr
 &= \cap\{ I(e) : I(e) \supseteq J \} \cr
 &= \ker\{ I\in \Omega : I \supseteq J \} \cr
 &=\ker(\hull(J)).
\endalign
$$
\qedd
\enddemo

It is interesting to see just what the topology on $\Omega$ is. Since $\Omega$
is in one-to-one
correspondence with the matrix units in system $\Cal E$, we can simply
transport the topology to
$\Cal E$ The closure of a point $e \in \Cal E$ is just the set
$$\overline{e} = \{ f\in{\Cal E} : f\leq_p e \},
$$
that is, the set of units smaller than $e$ in the partial order $\leq_p$. In
particular, this
topology is generally not $T_1$. For an arbitrary subset $F$ of $\Cal E$, we
can verify
$$\overline{F} = \{ f\in{\Cal E} : f\leq_p e, \hbox{ some } e \in F \}.
$$
This is a discrete topology, with a wedge condition on the closures.

\vskip .3cm

We now consider a certain subset of the almost finite dimensional triangular
algebras.
Namely,  let $\A$ be a norm-closed, strongly maximal triangular subalgebra of
an AF C*-algebra
$\B$. We will attempt to construct a topological space of meet-irreducible (and
nest-primitive) ideals which describe the ideal structure of $\A$, based on the
finite dimensional
constructions above.

Following Power in [16], fix an increasing sequence $\B_0\subseteq
\B_1\subseteq \B_2\subseteq
\ldots$ of finite dimensional C*-subalgebras of $\B$ with dense union, where
$\A_k = \B_k \cap \A$
is maximum triangular in $B_k$. The diagonal $\C = \A\cap\A^*$ is a canonical
masa, and is the
closure of the union of finite dimensional masas $\C_k = \A_k\cap\A_k^*$ in
$\B_k$. We choose a
system of matrix units
$\{ e_{ij}^k\}$ for each C*-algebra
$\B_k$ so that the matrix units in $\B_k$ are sums of matrix units in
$\B_{k+1}$, the
self-adjoint matrix units are precisely the ones in  the diagonal $\C_k$, and
each algebra $\A_k$
is a direct sum of  upper triangular matrices in $\B_k$.

We introduce the notion of an infinite chain of matrix units
$$e^{k_o} \to e^{k_o+1} \to e^{k_o+2} \to \ldots,$$
where each $e^{k}$ is a matrix unit in $\B_{k}$ from the system of matrix units
$\{ e_{ij}^k\}$,
 and the notation $e^k \to e^{k+1}$ indicates that $e^{k+1}$ is a summand of
the
matrix unit $e^k$. Since the system $\{ e_{ij}^k\}$  was constructed so that
each matrix unit in
$\B_k$ is a sum of matrix units in $\B_{k+1}$, a finite chain
$$e^{k_o} \to e^{k_o+1} \to e^{k_o+2} \to \ldots \to e^{k_o + n}$$
may be extended to an infinite chain by chasing down
 summands in the matrix unit system. With a finite choice of summands at each
stage $k$, there
are typically  uncountably many different ways of extending a chain.

If a matrix unit $e^k$ is in the maximal triangular algebra $\A$, then so are
all its summands.
Consequently, if  the first matrix unit in an infinite chain of matrix unit
lies in $\A$, then so
do all the units in the chain.

Let
$$e^{k_o} \to e^{k_o+1} \to e^{k_o+2} \to \ldots$$
be a chain of matrix units in $\A$, and for each $k\geq k_o$, let $I_k$ be the
largest ideal in
$\A_k$ not containing $e^k$. By the last section, $I_k$ is meet-irreducible and
(K4) in
$\A_k$. Observe that since $I_{k+1}$ does not contain $e^{k+1}$, the
intersection $I_{k+1}\cap
\A_k$ cannot contain $e^k$; since $I_k$ is the largest of such ideals, we have
$$I_k \supseteq I_{k+1}\cap \A_k.\tag"(*)"$$

In general, we don't have equality of these two ideals in $\A_k$ (an example is
given at the end
of this section). However, if
$I_k = I_{k+1}\cap
\A_k$ for all $k\geq k_o$, then by Theorem 2.5 of [12]  the closure of the
union of $I_k$  defines
an ideal $I = \overline{\cup I_k}$ in standard form.

\proclaim{Proposition 3.4} Let $\A = \overline{\cup \A_k}$ be a closed,
strongly maximal AF algebra
 and $\{ I_k \}$ a sequence of ideals constructed as above.
If $I_k = I_{k+1}\cap \A_k$, for all $k\geq k_o$, then the ideal $I =
\overline{\cup I_k}$  is both
(K4) and meet-irreducible.
\endproclaim

\demo{Proof}
To see that $I$ is (K4), let $J$ and $K$ be closed ideals in $\A$ with $I$
containing the
intersection $J\cap K$. Thus
$$I_k = I\cap\A_k \supseteq (J\cap K)\cap \A_k = (J\cap\A_k)\cap(K\cap\A_k) =
J_k \cap K_k.$$
Since $I_k$ is (K4) when $k\geq k_o$, we conclude that for all large $k$, $I_k$
contains either
$J_k$ or
$K_k$; in particular, either $I_k \supseteq J_k$ for infinitely many $k$, or
$I_k \supseteq K_k$
for infinitely many
$k$. Without loss of generality, say $I_k \supseteq J_k$ infinitely often.
Notice if $k'\leq k$,
where
$k$ is one of the indices where containment holds, then
$I_{k'} = I_k \cap \A_{k'} \supseteq J_k \cap \A_{k'} = J_{k'}$, so we conclude
$I_k \supseteq
J_k$ for all $k$, hence $I \supseteq J$. Thus $I$ is (K4).

Now by Proposition~1.3,  $I$ is  meet-irreducible as well.
\qedd
\enddemo

Note that since $I = \overline{\cup I_k}$ is in canonical form, we have that
$I_k = I\cap \A_k$ for
all $k\geq k_o$. In particular, since the matrix unit $e^k$ is not an element
of $I_k$, it
follows that every matrix unit in the chain is not in the ideal $I$.  This
provides a very rich
family of meet-irreducible ideals, provided we can assert the equality $I_k =
I_{k+1}\cap \A_k$
for all $k$. It is easy to verify that for both the standard embedding and
refinement embeddings
of upper triangular matrix algebras, we always have equality.
For if a matrix unit $f$ is not in $I_{k+1}\cap \T_{n_k}$, then writing $ f =
\sum f_i$ as a sum of
matrix units in $\T_{n_{k+1}}$, we find some $f_i$ is not in $I_{k+1}$, thus
$f_i\leq_p
e_{k+1}$. By the block structure of the standard and refinement embeddings, we
can
conclude $f \leq_p e_k$ and thus $f$ is not in $I_k$. Since ideals are sums of
the matrix
units they contain, we obtain the reverse containment of (*) above, namely
$I_k\supseteq
I_{k+1}\cap A_k$, and hence equality.

More generally, given an embedding
$\T_{n_1}\oplus\cdots\oplus\T_{n_k}
\rightarrow
\T_{m_1}\oplus\cdots\oplus\T_{m_{k'}}$, if each component map
$\T_{n_j}\rightarrow\T_{m_{j'}}$ is
either a standard or refinement embedding, or the zero map, then the  result
$I_k =
I_{k+1}\cap \A_k$ is similarly obtained. Thus in these cases, for any chain of
matrix units, we can
construct a corresponding meet-irreducible ideal. We note in the following
proposition that this
is a very rich family of meet-irreducibles.

\proclaim{Proposition 3.5} Let $\A$ be a strongly triangular AF algebra, where
each
non-zero component map of the embeddings
$\A_k
\to \A_{k+1}$ is either a standard or a refinement embedding. Then every closed
ideal in $\A$
equals the intersection of the meet-irreducible ideals which contain it.
Moreover, it suffices to
intersect only meet-irreducible ideals of the form constructed in
Proposition~3.4.
\endproclaim

\demo{Proof}
Let $J$ be a closed ideal in $\A$ and let $M$ be the intersection of all
meet-irreducible ideals
of the constructed form above, which contain $J$. Clearly $J$ is a subset of
$M$.

To show the reverse containment, suppose
$e^{k_o}$ is a matrix unit from $\A_{k_o}$ which is not in $J$. Since $e^{k_o}$
 is a
sum of matrix units in $\A_{k_o+1}$, there must be at least one summand
$e^{k_o+1}$ of $e^{k_o}$
in $\A_{k_o}$ which is not in the ideal $J$. Thus we may construct inductively
an infinite chain
of matrix units
$$e^{k_o} \to e^{k_o+1} \to e^{k_o+2} \to \ldots,$$
none of which lies in $J$. With $I_k$ the largest ideal in $\A_k$ not
containing $e_k$, we note
that $I_k \supseteq J\cap \A_k$, for all $k\geq k_o$. Thus the meet-irreducible
ideal $I =
\overline{\cup I_k}$ is of the constructed from above, $I$ contains $J$, and
$e^{k_o}$ is not in
$I$. Thus $e^{k_o}$ is not in the intersection of all such meet-irreducibles,
so $e^{k_o}$ is not
in $M$.

Since ideals equal the span of the matrix units they contain, we conclude $J =
M$.
\qedd
\enddemo

\proclaim{Corollary 3.6} With $\A$ as in Proposition 3.5, let $\Omega$ be the
set of all
meet-irreducible ideals constructed from unit chains $e^k \rightarrow e^{k+1}
\rightarrow \ldots$,
then the hull-kernel operation defines a topology on $\Omega$, and there is a
one-to-one
correspondence between closed sets in $\Omega$ and ideals in $\A$, via the hull
and kernel maps.
\endproclaim

Notice from the comments following Proposition~1.3, we can conclude that the
lattice of ideals in
$\A$ is distributive, and so the condition (K4) is equivalent to
meet-irreducible. This raises a
number of questions, including what does the space
$\Omega$ look like, and whether every meet-irreducible ideal in $\A$
corresponds to a point in
$\Omega$; that is,  is it of the special form constructed above. We will leave
these questions to
future work.

However, there is an interesting connection with nest representations which we
can describe.
With $X$ denoting the Gelfand space of the diagonal $\C$, we can build a
representation of $A$ on
the Hilbert space $l^2(X)$ by lifting the natural representation of each finite
algebra $\A_k$ on
its diagonal in the obvious way. To be explicit, we follow Power's description
in [16] of the
Gelfand space  and identify each point $x$ in $X$ with a sequence
$$ q^1_x \geq q^2_x \geq q^3_x \geq \ldots
$$
of diagonal projections $q_x^k$ in $\C_k$, where $q_x^k$ is the unique matrix
unit in $\C_k$ with
$x(q_x^k) = 1$. (The order $\geq$ is the usual range containment order for
commuting projections.)
Of course, given only the tail of a sequence
$$ q^k_x \geq q^{k+1}_x \geq q^{k+2}_x \geq \ldots
$$
the order condition $q_x^{k-1} \geq q_x^k$ determines uniquely the beginning of
the sequence
$$ q_x^1 \geq q_x^2 \geq \ldots \geq q^{k-1}_x   \geq q^k_x   \geq \ldots,
$$
so these sequences of diagonal projections corresponding to points in $X$ are
uniquely determined
by their tails. To build the representation of $\A$ on $l^2(X)$, let $\delta_x$
denote the basis
function in
$l^2(X)$ supported at the point $x$ in $X$, corresponding to some sequence of
units $q_x^1\geq
q_x^2\geq \ldots$, and for any matrix unit $f$ in $\A_k$, let $fx$ be the point
in $X$
corresponding to the tail sequence
$$f q_x^k f^* \geq f q_x^{k+1} f^* \geq f q_x^{k+2} f^* \geq \ldots,
$$
when this sequence doesn't terminate in zero. That is,  there is the
possibility that the sequence
collapses to
$f q_x^{k+n} f^* = 0$ for $n$ sufficiently large, in which case this sequence
does not represent
a point in $X$, so $fx$ is left undefined.

Define the map $\pi:\A \rightarrow \B(l^2(X))$ on matrix units by its action on
basis elements, as
$$\pi(f)\delta_x = \cases \delta_{fx}, &\text{if $fx \in X$ is defined }\\
                  0, &\text{if $f q_x^{k+n} f^* = 0$ for some $n\geq
0$}.\endcases
$$
This map $\pi$, when restricted to the units of $\A_k$, is just the natural
representation of
$\A_k$ on its diagonal, so it is easy to see that the linear extension of $\pi$
to
$\A_k$ creates a coherent representation of the finite subalgebras and thus
extends to a faithful
representation of $\A$.

To obtain a nest representation, fix a chain of matrix units
$$e^{k_o} \rightarrow e^{k_o + 1} \rightarrow e^{k_o+2} \rightarrow \ldots,
$$
 let $I_{k_o}, I_{k_o +1},\ldots$ be the corresponding sequence of
meet-irreducible ideals in
algebras
$\A_{k_o}, \A_{k_o+1},\ldots$, and let $X' \subseteq X$ be the subset of points
$x$ corresponding
to sequences of units
$$q_x^{k_o} \geq q_x^{k_o + 1} \geq q_x^{k_o + 2} \geq \ldots
$$
such that
$$q_x^k \not\in I_k, \qquad \hbox{ for all $k\geq k_o$}.
$$
Let $\pi'$ be the compression of the representation $\pi$ to the subspace
$l^2(X') \subseteq
l^2(X)$.

\proclaim{Proposition 3.7} Under the assumptions of Proposition 3.5, the linear
map $\pi'$ is a
nest representation of algebra $\A$, whose kernel is the meet-irreducible ideal
$I$ corresponding
to the chain of matrix units
$$e^{k_o} \rightarrow e^{k_o + 1} \rightarrow e^{k_o+2} \rightarrow \ldots.
$$
\endproclaim
\demo{Proof}
As in the construction of $\pi$, the map $\pi'$ restricted to any finite
subalgebra $\A_k$ is
just the compression of the natural representation of a triangular algebra
$\T_{n_1}\oplus\ldots\oplus
\T_{n_r}$ onto an interval of its diagonal, and is thus a (nest) representation
of $\A_k$. Taking
limits, we conclude $\pi'$ is a bounded representation of the AF algebra $\A$.
Note that for any
matrix unit $f$ in $\A_k$, we have $f\not\in \ker(\pi')$ if and only if
$f\leq_p e^k$, which is
equivalent to $f \not\in I_k = I\cap\A_k$. Since the ideal $I$ is determined by
the matrix units it
contains, we have $\ker(\pi') = I$.

To see that $\pi'$ is a nest representation, observe first that the
Peters-Poon-Wagner order on the
diagonal matrix units introduces a total order on the $X'$ via the relation
$$x \prec y \qquad \Leftrightarrow \quad
\cases q_x^k = q_y^k, &k = 1,2,\ldots n-1,\\
       q_x^n \prec q_y^n, &\text{some $n$}.\endcases
$$
Transitivity of this order relation follows from the observation that a
refinement embedding
component of $\A_k\rightarrow \A_{k+1}$ will preserve the Peters-Poon-Wagner
order for those
diagonal elements not in
$I_k$, while a standard embedding component is a one-to-one map on the diagonal
elements not in
$I_k$. Thus we may observe that $q_x^n \preceq q_y^n$ implies $q_x^k \preceq
q_y^k$ for all
$k\geq n$, from which transitivity follows. This order is total, since the
diagonal elements
$q_x^k, q_y^k$ in $\A_k$ which are not in $I$ form an interval along a diagonal
in $\A_k$, and
are thus totally ordered by $\preceq$.

The $\pi'$-invariant subspaces of $l^2(X')$ include the subspaces
$l^2[-\infty,x]$ of functions
supported on the interval $[-\infty,x] =\{ y\in X' : y\preceq x \}$, and more
generally, they are
the subspaces $l^2[-\infty,c]$, where $c$ represents a Dedekind cut in $X'$.
These form a nest via
the order
$\preceq $ on $X'$, thus $\pi'$ is a nest representation.
\qedd
\enddemo

Thus once again, we find a topological space of meet-irreducibles which
coincidently are
also nest-primitives, that completely describes the ideal structure of the
ideal space.

It is worth noting that if one considers more general *-extendable embeddings
of upper
triangular matrix algebras, the above construction need not work. In
particular, given a map
$\phi:\A_k \to \A_{k+1}$ and a meet-irreducible ideal $I_k\subseteq \A_k$, it
is not always
possible to choose a corresponding meet-irreducible ideal $I_{k+1} \subseteq
\A_{k+1}$ with
$I_k = I_{k+1}
\cap
\A_k$. For instance, consider the following embedding of $T_4$ into $T_8$,
$$
\left(
\matrix
a & b & c & d \\
  & e & f & g \\
  &   & h & i \\
  &   &   & j
\endmatrix
\right)
\to
\left(
\matrix
a & b &   &   & c & d &   &   \\
  & e &   &   & f & g &   &   \\
  &   & a & b &   &   & c & d \\
  &   &   & e &   &   & f & g \\
  &   &   &   & h & i &   &   \\
  &   &   &   &   & j &   &   \\
  &   &   &   &   &   & h & i \\
  &   &   &   &   &   &   & j \\
\endmatrix
\right),
$$
which is the refinement embedding, amplified by two. Let
$I_4$ be the meet-irreducible ideal in $T_4$ with corner $f$; that is, the
ideal with entries
$e,f,h$ set to zero. There are two possible choices for meet-irreducible ideals
$I_8\subseteq
\T_8$ with corner $f$. The first one (choosing the upper $f$ in $\T_8$) has
$I_8\cap\T_4$ equal to
the ideal with entries $a,b,e,f,h$ set to zero, which is strictly smaller than
$I_4$. The second
choice for $I_8$ has $I_8\cap\T_4$ equal to the ideal
with entries $e,f,h,i,j$ set to zero, which is also smaller than $I_4$.
By modifying this example with a twist in the last two pairs of rows/columns,
we obtain an
example where one choice of $I_8$ gives $I_4 = I_8\cap\T_4$ while the other
give $I_4 \neq
I_8\cap\T_4 = 0$.  We conclude that the general case of maximal triangular AF
algebras may be
quite complicated.

\head  Acknowledgements
\endhead
This research was supported in part by an
NSERC Individual Research Grant.

\enddocument